\documentclass[notitlepage]{amsart}

\usepackage{amsmath,amssymb,graphicx}
\usepackage{amsthm}
\usepackage[a4paper, top=3cm, bottom=2cm, left=3cm, right=3cm, marginparwidth=1.75cm]{geometry}
\usepackage[toc,page]{appendix}
\usepackage{enumitem}
\setlist[enumerate]{label*=\arabic*.}
\usepackage{braket}
\usepackage{parskip}

\newcommand{\norm}[1]{\left\lVert#1\right\rVert}
\newcommand{\abs}[1]{\left\lvert#1\right\rvert}

\newtheorem{theorem}{Theorem}[section]
\newtheorem{lemma}[theorem]{Lemma}
\newtheorem{corollary}[theorem]{Corollary}
\newtheorem{proposition}[theorem]{Proposition}

\theoremstyle{definition}
\newtheorem{definition}[theorem]{Definition}

\theoremstyle{remark}

\begin{document}

\title{Truncated Geometry on the Circle}
\author{E.~Hekkelman
}

\begin{abstract}
    In this letter we prove that the pure state space on the $n \times n$ complex Toeplitz matrices converges in the Gromov-Hausdorff sense to the state space on $C(S^1)$ as $n$ grows to infinity, if we equip these sets with the metrics defined by the Connes distance formula for their respective natural Dirac operators. A direct consequence of this fact is that the set of measures on $S^1$ with density functions $c \prod_{j=1}^n (1-\cos(t-\theta_j))$ is dense in the set of all positive Borel measures on $S^1$ in the weak$^*$ topology.
\end{abstract}

\maketitle

\section{Introduction}
In the framework of noncommutative geometry~\cite{Connes1994}, we study the pure states on the Toeplitz matrices as a metric space. The results of this letter can be understood with little prior knowledge, however. With $n\times n$ Toeplitz matrices we mean those matrices in $M_n (\mathbb{C})$ for which each descending diagonal is constant, i.e. the matrices $T$ for which $T_{i,j} = T_{i+1, j+1}$. We denote this set by $C(S^1)^{(n)}$, as these Toeplitz matrices can also be seen as a truncation of $C(S^1)$ which we will explain later on. The reader need only know that these matrices inherit a structure of positivity from $M_n(\mathbb{C})$ which allows us to talk about positive linear functionals. Likewise, states on the Toeplitz system can be defined as those positive linear functionals $\varphi$ for which $\varphi(I)=1$, and pure states as the extreme points in the state space. We denote these spaces as $\mathcal{S}(C(S^1)^{(n)})$ and $\mathcal{P}(C(S^1)^{(n)})$, respectively. The states on $C(S^1)^{(n)}$ can be equipped with a metric that metrises the weak$^*$-topology via the Connes distance formula which is defined as \[d_n(\varphi, \psi) = \sup_{T \in C(S^1)^{(n)}} \{\abs{\varphi(T)-\psi(T)} : \norm{[D_n,T]} \leq 1 \},\] where $D_n$ denotes the matrix $D_n := \operatorname{diag}(1, \dots, n)$. Importantly, this is completely analogous to a formula that defines the Monge-Kantorovich metric on the state space of $C(S^1)$, \[d(\varphi, \psi) = \sup_{f \in C^\infty(S^1)} \{\abs{\varphi(f)-\psi(f)} : \lVert[D,f]\rVert \leq 1 \},\] where $D= -i\frac{d}{dx}$ is the Dirac operator on $S^1$, which metrises the weak$^*$-topology on  $\mathcal{S}(C(S^1))$~\cite{rieffel1999metrics}. In this language, the main result of this letter is that the metric spaces $(\mathcal{P}(C(S^1)^{(n)}), d_n)$ converge in the Gromov-Hausdorff sense to the metric space $(\mathcal{S}(C(S^1)), d)$.

The reader who is familiar with noncommutative geometry may be interested to know how the above relates to prior work in that area. Related to physical applications of noncommutative geometry, extensive study has been made of spectral geometry for which only part of the spectral data is available~\cite{ConnesSuijlekom, DAndrea, vansuijlekom2020gromovhausdorff}. Indeed, detectors in physical experiments are limited and thus only give information up to a certain energy level with finite resolution. Suggested first by F. D’Andrea, F. Lizzi, and P. Martinetti~\cite{DAndrea} and put in the language of operator systems in an article by A. Connes and W.D. van Suijlekom~\cite{ConnesSuijlekom}, the correct theoretical framework for such truncations seems to be to truncate a spectral triple $(\mathcal{A}, H, D)$ with a projection $P$ onto a part of the spectrum of $D$ resulting in the triple $(P\mathcal{A}P, PH, PD)$. The Toeplitz matrices arise as exactly such a truncation of the algebra $\mathcal{A}=C(S^1)$ in the spectral triple of the circle. Note also that a first result on Gromov-Hausdorff convergence of state spaces in this context is put forward in~\cite{DAndrea}. A more detailed explanation of the connection between the results presented in this letter and noncommutative geometry can be found in the Master's thesis~\cite{Scriptie}.

Before continuing on, the author would like to thank Walter van Suijlekom for his invaluable guidance.

\section{Preliminaries}
These preliminaries consist of two parts, namely preliminaries on Gromov-Hausdorff convergence and on truncated geometry on the circle.

\subsection{Gromov-Hausdoff convergence}
\label{GH}
In the spirit of work by M. Rieffel~\cite{Rieffel20041}, we are interested in studying limits of sequences of metric spaces in the Gromov-Hausdorff sense. The Gromov-Hausdorff distance is a pseudo-metric on the class of all compact metric spaces, and is zero if and only if two spaces are isometric~\cite[Theorem~7.3.30]{burago}. Hence it functions a metric on the isomorphism classes of compact metric spaces, although one has to be careful to dodge set-theoretic paradoxes in such a description. In any case, it gives a useful notion of convergence and there are several techniques that can be used to prove Gromov-Hausdorff convergence of sequences of compact metric spaces.

For a subset $S$ in a metric space, denote the $r$-neighbourhood of $S$ by $U_r(S)$, i.e. \[U_r(S) = \bigcup_{x\in S} B_r(x),\] where $B_r(x)$ is the open ball of radius $r$ and center $x$.
\begin{definition}
    Let $A$ and $B$ be subsets of a metric space. The \textit{Hausdorff distance} between $A$ and $B$ is defined as \[d_H(A,B) = \inf\{ r > 0 : A \subseteq U_r(B) \text { and } B \subseteq U_r(A)\} .\]
\end{definition}

\begin{definition}
\label{GMMetric}
    Let $X$ and $Y$ be metric spaces. The \textit{Gromov-Hausdorff distance} between $X$ and $Y$ is defined as the infimum of all $r > 0$ such that there exists a metric space $Z$ with subsets $X', Y' \subseteq Z$ isometric to $X$ and $Y$ respectively with $d_H(X', Y') < r$, where $d_H(X', Y')$ is the Hausdorff distance between $X'$ and $Y'$. We will denote the Gromov-Hausdorff distance by $d_{GH}(X, Y)$ in the rest of this letter as well.
\end{definition}

These definitions are exactly~\cite[Definition~7.3.1]{burago} and~\cite[Definition~7.3.10]{burago}.

\begin{definition}\label{correspondence}
    Let $X$ and $Y$ be two sets. A \textit{total onto correspondence} between $X$ and $Y$ is a set $\mathfrak{R} \subseteq X \times Y$ such that for every $x \in X$ there exists at least one $y \in Y$ with $(x,y) \in \mathfrak{R}$ and similarly for every $y \in Y$ there exists an $x \in X$ with $(x,y) \in \mathfrak{R}$.
\end{definition}

Note that the above definition is simply called a correspondence in~\cite[Chapter~7]{burago}, although this is not usual terminology. To prevent any confusion we will stick to calling the above \textit{total onto} correspondences.

\begin{definition}
    Let $\mathfrak{R}$ be a total onto correspondence between metric spaces $X$ and $Y$. The \textit{distortion} of $\mathfrak{R}$ is defined by \[\operatorname{dis } \mathfrak{R} = \sup \{\abs{d_X(x, x') - d_Y(y, y')} : (x,y), (x', y') \in \mathfrak{R} \}.\]
\end{definition}

\begin{theorem}
\label{GHCorrespondence}
For any two metric spaces $X$ and $Y$ \[d_{GH}(X,Y) = \frac{1}{2} \inf_{\mathfrak{R}}(\operatorname{dis } \mathfrak{R}),\] where the infimum is taken over all total onto correspondences $\mathfrak{R}$ between $X$ and $Y$.
\begin{proof}
See~\cite[Theorem~7.3.25]{burago}.
\end{proof}
\end{theorem}

As a direct corollary, there is a similar approach of using $\varepsilon$-isometries.

\begin{definition}
Let $X$ be a metric space and $\varepsilon > 0$. A set $S \subseteq X$ is called an \textit{$\varepsilon$-net} if $\operatorname{dist}(x,S) \leq \varepsilon$ for every $x \in X$.
\end{definition}

\begin{definition}{Let $X$ and $Y$ be metric spaces and $f: X \rightarrow Y$ an arbitrary map. The \textit{distortion} of $f$ is defined by \[\operatorname{dis}f = \sup_{x_1, x_2 \in X} |d_Y (f(x_1), f(x_2)) - d_X (x_1, x_2)|.\]}
\end{definition}

\begin{definition}
Let $X$ and $Y$ be metric spaces and $\varepsilon > 0$. A (possibly non-continuous) map $f: X \rightarrow Y$ is called an \textit{$\varepsilon$-isometry} if $\operatorname{dis}f \leq \varepsilon$ and $f(X)$ is an $\varepsilon$-net in $Y$.
\end{definition}

\begin{corollary}{\label{epsilonisometry}Let $X$ and $Y$ be metric spaces and $\varepsilon > 0$. \begin{enumerate}
    \item If $d_{GH}(X, Y) < \varepsilon$, then there exists a $2\varepsilon$-isometry from $X$ to $Y$.
    \item If there exists an $\varepsilon$-isometry from $X$ to $Y$, then $d_{GH}(X, Y) < 2\varepsilon$.
\end{enumerate}}
\begin{proof}
See~\cite[Corollary~7.3.28]{burago}.
\end{proof}
\end{corollary}

More details on the Gromov-Hausdorff metric can be found in Chapter 7 of \textit{A Course in Metric Geometry}~\cite{burago} by D. Burago, I. Burago and S. Ivanov.

\subsection{Truncated Geometry on the Circle}
\label{TruncGeomCirc}
As mentioned in the introduction, the Toeplitz matrices arise as a truncation of the algebra of continuous functions $C(S^1)$. Indeed, $C(S^1)$ has a basis $\{e_n(t) = e^{int} \}_{n \in \mathbb{Z}}$, consisting of all eigenfunctions of the Dirac operator $D = -i \frac{d}{dx}$ on $S^1$. We can then consider the projection onto $\operatorname{span}_\mathbb{C}\{e_1, \cdots, e_n\}$, denote this projection by $P_n$. For any function $f \in C^\infty(S^1)$, the action of $P_nfP_n$ on the finite dimensional Hilbert space $\mathrm{span}_\mathbb{C} \{e_k\}_{k=1}^n$ can be represented as the $n \times n$ Toeplitz matrix \[\left(\begin{array}{ccccc}
     \hat{f}(0)& \hat{f}(-1) & \hat{f}(-2) & \cdots & \hat{f}(-n+1) \\
     \hat{f}(1) & \hat{f}(0) & \hat{f}(-1) & \cdots & \hat{f}(-n+2)\\ 
    \hat{f}(2) & \hat{f}(1) & \hat{f}(0) & \cdots & \hat{f}(-n+3)\\ 
    \vdots & \vdots & \vdots & \ddots & \vdots \\
    \hat{f}(n-1) & \hat{f}(n-2) & \hat{f}(n-3) & \cdots & \hat{f}(0)\\ 
\end{array} \right).\]
This is the reason why we denote the set of Toeplitz matrices by $C(S^1)^{(n)} $.

It has already been proven in~\cite{vansuijlekom2020gromovhausdorff} that the state spaces $\mathcal{S}(C(S^1)^{(n)})$ converge to $\mathcal{S}(C(S^1))$ in the Gromov-Hausdorff sense. To prove our result, a main ingredient of~\cite{vansuijlekom2020gromovhausdorff} will be used in this letter as well, which is the map \begin{align*}
    R_n: C(S^1) & \rightarrow C(S^1)^{(n)}\\
    f \quad & \mapsto P_n f P_n.
\end{align*}

\begin{lemma}\label{R_n}
The map \begin{align*}
    R_n^*: \mathcal{S}(C(S^1)^{(n)}) &\rightarrow \mathcal{S}(C(S^1))\\
    \tau \quad \quad &\mapsto \quad \tau \circ R_n
\end{align*}
is well-defined and satisfies $\operatorname{dis} R_n^* \rightarrow 0$ as $n \rightarrow \infty$.
\begin{proof}
See~\cite{vansuijlekom2020gromovhausdorff}.
\end{proof}
\end{lemma}

Next, let us give a description of the pure states on $C(S^1)^{(n)}$. This has already been done via a duality of the Toeplitz operator system with the Fej\'er-Riesz operator system~\cite{ConnesSuijlekom}, but here we will use a more direct approach.

A very useful ingredient for this is the following decomposition theorem dating from 1911 proven by C. Caratheodory and L. Fej\'er~\cite{caratheodory}. We introduce the notation \[f_z = \frac{1}{\sqrt{n}} \ \big(\begin{array}{ccccc}
     1 & z & z^2 & \cdots & z^{n-1}
\end{array}\big) \in \mathbb{C}^n,\] which is a column of a so-called Vandermonde matrix.

\begin{theorem}\label{vandermonde}
Any positive Toeplitz matrix $T \in C(S^1)^{(n)}$ of rank $r \leq n-1$ can be uniquely decomposed as $T = \sum_{k=1}^{r} d_k \ket{f_{\lambda_k}} \bra{f_{\lambda_k}}$ where $d_1, \dots, d_r > 0$ and $\lambda_1, \dots, \lambda_r \in S^1$. This is called the \textit{Vandermonde decomposition}. If the rank of $T$ is $n$, this decomposition is still possible but not unique.
\begin{proof}
See~\cite{vandermonde}.
\end{proof}
\end{theorem}

While we use this classical theorem to classify the pure states on the Toeplitz operator system, the same result can also be derived by the aforementioned operator system duality, see~\cite[Theorem~4.14]{ConnesSuijlekom}.

\begin{proposition}\label{purestates1}
A state on $C(S^1)^{(n)}$ is pure if and only if it is a vector state $\varphi_\xi: T \mapsto \langle \xi, T\xi \rangle,$ where $\xi = (\xi_0, \xi_1, \dots, \xi_{n-1}) \in \mathbb{C}^n$ is a unit vector such that the polynomial $Q_\xi(z) := \sum_k \xi_k z^{n-k-1}$ has all its zeroes on $S^1$.\begin{proof}
First of all, any pure state on $C(S^1)^{(n)}$ admits an extension to a pure state on $M_n(\mathbb{C})$~\cite[Fact~2.9]{ConnesSuijlekom}. Since all pure states on $M_n(\mathbb{C})$ are vector states, any pure state on $C(S^1)^{(n)}$ must likewise be of the form $\varphi_\xi: T \mapsto \langle \xi, T\xi \rangle$ for some unit vector $\xi \in \mathbb{C}^n$.

By the Vandermonde decomposition of Toeplitz matrices, we can write any positive $T \in C(S^1)^{(n)}$ in the form \[\sum_{k=1}^{r} d_k \ket{f_{\lambda_k}} \bra{f_{\lambda_k}}\] with $d_1, \dots, d_{r} \geq 0$, this is Theorem~\ref{vandermonde}. Observe that on $S^1$, we have $\overline{\lambda} = \lambda^{-1}$ and so \[\langle f_\lambda, \xi \rangle = \sum_{k=0}^{n-1} \xi_k \lambda^{-k} = \lambda^{-n+1} Q_\xi(\lambda).\] Therefore, \[\varphi_\xi(T) = \langle \xi, T\xi \rangle = \sum_{k=1}^{r} d_k \abs{\langle f_{\lambda_k}, \xi\rangle}^2 = \sum_{k=1}^{r} d_k \abs{Q_\xi(\lambda_k)}^2.\] 

Now let us investigate for what $\xi$ this state is pure. Note that $\varphi_\xi$ is pure if and only if \[\varphi_\omega(T) \leq \varphi_\xi(T) \quad \forall 0\leq T\in C(S^1)^{(n)}\] implies that $\omega \in \mathbb{C}\xi$ (see~\cite[p.~144]{Murphy}), and by the above calculation this is equivalent to \[\abs{Q_\omega(z)} \leq \abs{Q_\xi(z)} \quad \forall z \in S^1\] implying that $\omega \in \mathbb{C}\xi$.

We claim that if $\abs{Q_\omega(z)} \leq \abs{Q_\xi(z)}$ for all $z \in S^1$ and $\lambda \in S^1$ is a root of $Q_\xi$, then $\lambda$ must also be a root of $Q_\omega$ with at least the same multiplicity. The first part of this claim is immediate. To prove the second part, suppose the multiplicity of the root $\lambda$ is $m$ for $Q_\omega$ but strictly more than that for $Q_\xi$. Then we must still have \[\abs{\frac{Q_\omega(z)}{(z-\lambda)^m}} \leq \abs{\frac{Q_\xi(z)}{(z-\lambda)^m}}\] on $S^1$, but with the right-hand side evaluating to zero at $z = \lambda$ whereas the left-hand side does not, we arrive at a contradiction. 

Therefore, if $Q_\xi$ has all its roots on $S^1$, the above observation combined with the fact that $Q_\omega$ and $Q_\xi$ are polynomials of the same degree leads to the conclusion that $\abs{Q_\omega(z)} \leq \abs{Q_\xi(z)} \forall z\in S^1$ implies that $Q_\omega \in \mathbb{C}Q_\xi$. Hence indeed $\varphi_\xi$ is pure, which completes one direction of the proposition.

For the other direction, suppose that $\xi \in \mathbb{C}^n$ is a vector such that the polynomial $Q_\xi$ has roots $\lambda_1, \cdots, \lambda_n $ (counted with multiplicities), of which say $\lambda_n$ is not an element of $S^1$. Then $\abs{z-\lambda_n}$ attains a minimum on $S^1$ which is strictly greater than zero, say $\delta > 0$. Choose any $\lambda \in S^1$ and note that $\abs{z-\lambda} \leq 2$ on $S^1$. Then \begin{align*}
    \abs{Q_\xi(z)}^2 &= \big \lvert c \prod_{k=1}^n (z-\lambda_k)\big\rvert^2\\
    & \geq \frac{\delta^2}{4}\big\vert c (z-\lambda) \prod_{k=1}^{n-1} (z-\lambda_k)\big\rvert^2,
\end{align*}
so the polynomial $\frac{\delta}{2} c (z-\lambda) \prod_{k=1}^{n-1} (z-\lambda_k) =: \sum_k \omega_k z^{n-k-1}$ corresponds to some $\omega \in \mathbb{C}^n$ with the property that \[\langle \omega, T\omega \rangle \leq \langle \xi, T\xi \rangle\] for all positive $T \in C(S^1)^{(n)}$. Clearly, $Q_\omega$ is not a scalar multiple of $Q_\xi$ so $\omega \not \in \mathbb{C}\xi$, and hence the vector state $\varphi_\xi$ is not pure.
\end{proof}
\end{proposition}

\begin{corollary}{\label{purestates2} The pure states of $C(S^1)^{(n)}$ are exactly the linear functionals of the form $P_nfP_n \mapsto \! \int_{S^1} f \ \abs{Q_\xi}^2 \ d\lambda$, where $\xi \in \mathbb{C}^n$ is a unit vector and $Q_\xi:= \sum_{k=0}^{n-1}\xi_k z^{n-k-1}$ has all its roots on $S^1$.}
\begin{proof}
According to Proposition~\ref{purestates1}, the pure states on $C(S^1)^{(n)}$ are given by $T \mapsto \langle \xi, T\xi \rangle,$ where the unit vector $\xi \in \mathbb{C}^{n}$ is such that the polynomial $Q_\xi$ has all its roots on $S^1$. A short calculation gives that \[\langle \xi, P_nfP_n \xi \rangle = \sum_{\abs{j}\leq n-1} (\xi^* * \xi)_j \ \hat{f}(-j),\]
where $(v * w)_j = \sum_{k=\max(0, -j)}^{\min(n-1,n-j-1)} v_{-k} w_{j+k}$ is the discrete convolution product and $(v^*)_j = \overline{v_{n-j-1}} $. By noting that $(\xi^* * \xi)_j$ is exactly the $j$th Fourier coefficient of $\abs{Q_\xi}^2$ we can use the Plancherel Theorem~\cite[Theorem~9.13]{PapaRudin} to conclude that \[\langle \xi, P_nfP_n \xi \rangle = \int_{S^1} f \ \abs{Q_\xi}^2 \ d\lambda.\qedhere\] 
\end{proof}
\end{corollary}

This corollary essentially characterises the pure states $\tau$ of $C(S^1)^{(n)}$ by the Radon-Nikodym derivative of $R_n^*\tau$ as a state (i.e. probability measure) on $C(S^1)$. This works since the Fourier basis $\{e_n\}_{n \in \mathbb{Z}}$ forms an orthonormal basis of $C(S^1)$, hence for a function $g \in \mathrm{span}_\mathbb{C} \{e_1, \dots e_n \}$ there is some ambivalence in considering a state $\tau: P_nfP_n \mapsto \int_{S^1} fg d\lambda$ on $C(S^1)^{(n)}$ or its pullback $R_n^*\tau: f \mapsto \int_{S^1} fg d\lambda$. We will exploit this, although responsibly in order to prevent confusion.

\textbf{Notation.} When considering a pure state $\tau$ on $C(S^1)^{(n)}$, we will somewhat abusively refer to the Radon-Nikodym derivative of $R_n^*\tau$, with respect to the normalised Haar-measure $d\lambda$ on $S^1$, as $\frac{d\tau}{d\lambda}$ (instead of $\frac{dR_n^*\tau}{d\lambda}$) because this function uniquely defines the pure state $\tau$. In the other way around, if $f$ is a function of the form such that it defines a pure state on $C(S^1)^{(n)}$, we will denote that pure state $\tau_{f}$. In summary, $\frac{d\tau_f}{d\lambda} = f$.

\begin{proposition}{\label{Pure state form}If $\tau$ is a pure state on $C(S^1)^{(n)}$, then \[\frac{d\tau}{d\lambda}(t) = \abs{Q_\xi}^2(e^{it}) = c \prod_{j=1}^{n-1} (2-2\cos(t-\theta_j)),\] where $e^{i\theta_j}$ are the roots of the polynomial $Q_\xi := \sum_k \xi_k z^{n-k-1}$ and $c\in\mathbb{R}$ is a scaling factor such that $\frac{d\tau}{d\lambda}$ integrates to 1. Likewise, any function of this form defines a pure state.}
\begin{proof}
According to Corollary~\ref{purestates2}, any pure state $\tau$ on $C(S^1)^{(n)}$ corresponds to a function of the form $\abs{Q_\xi}^2$. Denoting the roots of $Q_\xi$ by $e^{i\theta_j}$, there must be a constant $c$ such that \begin{align*}
    \abs{Q_\xi(e^{it})}^2   &= c\prod_{j=1}^{n-1} \abs{e^{it} - e^{i \theta_j}}^2\\
    &= c\prod_{j=1}^{n-1} (2 - 2 \cos(t-\theta_j)).
\end{align*}
Furthermore, $\abs{Q_\xi}^2$ must integrate to 1 because $\norm{Q_\xi}_2 = \norm{\xi}$ since the Fourier transform is unitary.

For the other way around, the above proves that any function of the form $ c\prod_{j=1}^{n-1} (2 - 2 \cos(t-\theta_j))$ is equal to $\abs{Q_\xi(e^{it})}^2$ for some vector $\xi \in \mathbb{C}^{n}$ such that $Q_\xi$ has all its zeroes on $S^1$. And again because $\norm{Q_\xi}_2 = \norm{\xi}$, $\abs{Q_\xi}^2$ integrating to 1 implies that $\xi$ is a unit vector. According to Corollary~\ref{purestates2}, this function therefore indeed defines a pure state on $C(S^1)^{(n)}$.
\end{proof}
\end{proposition}

\section{Convergence to $\mathcal{S}(C(S^1))$}
\label{ConvergenceSection}
We will now prove the convergence of $\mathcal{P}(C(S^1)^{(n)})$ to $\mathcal{S}(C(S^1))$ by establishing $\varepsilon$-isometries between these spaces, employing Corollary~\ref{epsilonisometry}. The candidate maps we propose are the maps $R_n^*$ for which we already have that $\operatorname{dis} R_n^* \rightarrow 0$ (Lemma~\ref{R_n}). If we can now establish that $R_n^* (\mathcal{P}(C(S^1)^{(n)}))$ is an $\varepsilon$-net if we choose $n$ large enough, Gromov-Hausdorff convergence follows directly. First, we will prove that all states on the circle can be approximated by the pullback of pure states in $\mathcal{P}(C(S^1)^{(n)})$. Next, we check that this can be done uniformly in $n$.

\subsection{Approximating states}
In this subsection, three steps of increasing difficulty will be carried out to prove that any state on the circle can be approximated by the pullbacks of pure states on the truncated circle. 

\begin{enumerate}
    \item We approximate any state $\psi$ on $C(S^1)$ by a state $\sum_{i=1}^m t_i \text{ev}_{\lambda_i}$ on $C(S^1)$, which is a convex combination of evaluations at the $m$-roots of unity;
    \item We then approximate any such convex combination of evaluations by a state \[\sum_{i=1}^m \frac{ \frac{d\varphi}{d\lambda}(\lambda_i)}{\sum_{j=1}^m \frac{d\varphi}{d\lambda}(\lambda_j) } \text{ev}_{\lambda_j}\] on $C(S^1)$, where $\lambda_j$ are the $m$-roots of unity and $\varphi$ is a pure state on $C(S^1)^{(n)}$;
    \item Finally, we approximate any state of that particular form by a state $R_n^* \chi$ on $C(S^1)$ where $\chi$ is a pure state on $C(S^1)^{(n)}$.
\end{enumerate}

Recall that a pure state $\varphi$ on $C(S^1)^{(n)}$ is uniquely characterised by the Radon-Nikodym derivative with respect to the normalised Haar-measure on $S^1$ of $R_n^*\varphi$, which is a state (i.e. a probability measure) on $C(S^1)$, and that we denote this Radon-Nikodym derivative $\frac{d\varphi}{d\lambda}$ instead of $\frac{dR_n^*\varphi}{d\lambda}$ to ease notation. See also Subsection~\ref{TruncGeomCirc}. For the third and most important step, the essential property of the pure state space we will exploit is that it is possible to multiply these Radon-Nikodym derivatives to construct new pure states, as can be seen from the form of these functions in Proposition~\ref{Pure state form}.

The first step in the outlined scheme is by far the easiest. Because the roots of unity are dense in $S^1$ and the standard topology on $S^1$ coincides with the weak$^*$-topology on the pure states, which in turn is induced by the Monge-Kantorovich metric~\cite{rieffel1999metrics}, it follows that the set \[\left\{ \sum_{i=1}^m t_i \text{ev}_{\lambda_i} \ : \ m \in \mathbb{N}, \ 0 \leq t_i \leq 1, \sum_{i=1}^m t_i = 1, \lambda_i^m = 1 \right\},\] i.e. convex combinations of evaluations at the roots of unity, is dense in $\mathcal{S}(C(S^1))$ with respect to the topology induced by the Monge-Kantorovich metric. 

The second step can be done in a single lemma. In spirit, this lemma should be compared with the assertion that one can always fit a polynomial $P$ of degree $n-1$ (or higher, if desired) through $n$ prescribed points. Here, we must choose a pure state $\varphi$ such that the function \[\frac{\frac{d\varphi}{d\lambda}}{\sum_{i=1}^m \frac{d\varphi}{d\lambda}(\lambda_i)}\] evaluates (approximately) to $t_i$ at the points $\lambda_i$.

\begin{lemma}{ \label{approx}Let $\lambda_1, ..., \lambda_m$ be the solutions of $\lambda^m = 1$ (the $m$-roots of unity), and take any state of the form $\sum_{i=1}^m t_i \text{ev}_{\lambda_i}$ with $\sum_{i=1}^m t_i = 1$ and $t_i \geq 0$ for all $i$. Then for every $l \in \mathbb{Z}_{\geq 0}$ and $\varepsilon >0$ we can find $\varphi \in \mathcal{P}(C(S^1)^{(m+1+l)})$ such that on $C(S^1)$ \[d\left(\sum_{i=1}^m t_i \text{ev}_{\lambda_i}, \sum_{j=1}^m \frac{\frac{d\varphi}{d\lambda}(\lambda_j)}{\sum_{i=1}^m \frac{d\varphi}{d\lambda}(\lambda_i)}   \text{ev}_{\lambda_j} \right) < \varepsilon.\]}\begin{proof}
Let us first prove the case $l = 0$. Take the pure state $\varphi_N$ on $C(S^1)^{m+1}$ defined by $\frac{d\varphi_N}{d\lambda}(t) = c\prod_{j=1}^m (1-\cos(t-\lambda_j + \sqrt{\frac{2t_j}{N}}))$. This indeed corresponds to a pure state in $\mathcal{P}(C(S^1)^{(m+1)})$ according to Proposition~\ref{Pure state form}. Consider what happens at the points $\lambda_i$ if we take $N$ large. 

By Taylor expansion of the factors, it is a quick calculation to see that
\begin{align*}
    \frac{d\varphi_N}{d\lambda}(\lambda_i) &= c\left(\frac{t_i}{N}+ \mathcal{O}\left(\frac{1}{N^{3/2}}\right)\right) \left( \prod_{j \not = i} \left((1-\cos(\lambda_i-\lambda_j)) + \mathcal{O}\left(\frac{1}{\sqrt{N}}\right)\right) \right)\\
    &= c\frac{t_i}{N}\prod_{j \not = i} (1-\cos(\lambda_i-\lambda_j)) + \mathcal{O}\left(\frac{1}{N^{3/2}}\right).
\end{align*}

Notice that $c\prod_{j \not = i} (1-\cos(\lambda_i-\lambda_j))$ has the same value for all $\lambda_i$ by symmetry. Hence if we pass these values to the projective space $\mathbb{R}P^{m-1}$ we end up with the ratio \begin{align*}
    \left[\frac{d\varphi_N}{d\lambda}(\lambda_1) : \cdots : \frac{d\varphi_N}{d\lambda}(\lambda_m)\right] &= \left[\frac{t_1}{N} + \mathcal{O}\left(\frac{1}{N^{3/2}}\right): \cdots : \frac{t_m}{N} + \mathcal{O}\left(\frac{1}{N^{3/2}}\right)\right]\\
    &=  \left[t_1 + \mathcal{O}\left(\frac{1}{\sqrt{N}}\right) : \cdots : t_m + \mathcal{O}\left(\frac{1}{\sqrt{N}}\right) \right].
\end{align*}

It then follows that the vectors \[\frac{1}{\sum_{j=1}^m \frac{d\varphi_N}{d\lambda}(\lambda_j)}\left(\frac{d\varphi_N}{d\lambda}(\lambda_1), \dots, \frac{d\varphi_N}{d\lambda}(\lambda_m)\right)\] converge to $(t_1, \dots, t_m)$ in $\mathbb{R}^m$ as $N \rightarrow \infty$. Therefore the states $\frac{1}{\sum_{j=1}^m \frac{d\varphi_N}{d\lambda}(\lambda_j)} \sum_{j=1}^m \frac{d\varphi_N}{d\lambda}(\lambda_j)   \text{ev}_{\lambda_j}$ converge to $\sum_{j=1}^m t_j \text{ev}_{\lambda_j}$ in the weak$^*$-topology. Again we can use that the Monge-Kantorovich metric induces the weak$^*$-topology to conclude that we can choose $N$ such that $$d\left(\sum_{i=1}^m t_i \text{ev}_{\lambda_i}, \sum_{j=1}^m \frac{\frac{d\varphi_N}{d\lambda}(\lambda_j)}{\sum_{i=1}^m \frac{d\varphi_N}{d\lambda}(\lambda_i)}   \text{ev}_{\lambda_j} \right) < \varepsilon.$$

The cases $l \geq 1$ follow more or less immediately. If we choose some point $\mu$ on the circle that is not equal to any of the $\lambda_j$, we can guarantee that $1-\cos(t-\mu)$ has no roots in the points $\lambda_1, \dots, \lambda_m$. For any $l \in \mathbb{N}$, we can take $\varphi_N \in \mathcal{P}(C(S^1)^{(m)})$ such that the ratio $ [\frac{d\varphi_N}{d\lambda}(\lambda_1) : \cdots : \frac{d\varphi_N}{d\lambda}(\lambda_m)]$ is arbitrarily close to \[\left[\frac{t_1}{(1-\cos(\lambda_1 - \mu))^l} : \cdots : \frac{t_m}{(1-\cos(\lambda_m - \mu))^l}\right]\] by the argument for the case $l=0$ above. Then $\frac{d\varphi_N}{d\lambda} (1-\cos(t- \mu))^l$ defines the pure state (up to scaling) in $\mathcal{P}(C(S^1)^{(m+1+l)})$ that satisfies the statement in the lemma. 
\end{proof}
\end{lemma}

For the third and final step, we need to prove that the states of this type can be approximated by the pullback of some pure state on the truncated circle. To accomplish that, we will need the following propositions.

\begin{proposition}{Let $K\subseteq X$ be some compact subset of $\mathbb{R}^n$ and let $f \in C(K)$ be a positive function attaining its maximum in the unique point $x_0$. Then the sequence of linear functionals $\left(\tau_n\right)_{n \in \mathbb{N}}$ defined by $\tau_n: g \mapsto \int_K \frac{f^n}{\norm{f^n}_1} g dx$, converges to $ \text{ev}_{x_0}$ in the weak$^*$-topology on $C(K)^*$.}
\label{Bumpfunctions}
\begin{proof}
Denote the maximum of $f$ by $M$. For every $\varepsilon > 0$, $f^{-1} (M - \varepsilon, M]$ is an open neighbourhood of $x_0$, denote this by $U_\varepsilon$. Outside this neighbourhood $\frac{f^n}{\norm{f^n}_1} \xrightarrow{n \rightarrow \infty} 0$ uniformly, since \[ \norm{f^n}_1 \geq \int_{U_{\varepsilon/2}} f^n dx > \abs{U_{\varepsilon/2}} (M - \varepsilon/2)^n, \] and so for $x \not \in U_\varepsilon$ \[ \frac{f^n}{\norm{f^n}_1} (x) \leq \frac{1}{\abs{U_{\varepsilon/2}}} \left(\frac{M-\varepsilon}{M - \varepsilon/2}\right)^n.\]
Therefore, 

\begin{align*}
    \abs{\tau_n(g) - \text{ev}_{x_0}(g)}&=\abs{\int_K \frac{f^n}{\norm{f^n}_1} g dx - \text{ev}_{x_0}(g)}\\
    &= \abs{\int_K \frac{f^n}{\norm{f^n}_1} (g - g(x_0)) dx}\\
    & \leq \int_{K-U_\varepsilon} \frac{f^n}{\norm{f^n}_1} \abs{g - g(x_0)}dx + \abs{\int_{U_\varepsilon} \frac{f^n}{\norm{f^n}_1} (g-g(x_0)) dx}\\
    &\leq \underbrace{\int_{K-U_\varepsilon} \frac{f^n}{\norm{f^n}_1} \abs{g - g(x_0)}dx}_{\xrightarrow{n \rightarrow \infty} 0} + \underbrace{\sup_{x \in U_\varepsilon} \abs{g(x)-g(x_0)}}_{ \xrightarrow{\varepsilon \rightarrow 0} 0} \underbrace{\int_{U_\varepsilon} \frac{f^n}{\norm{f^n}_1}dx}_{\leq 1}.
\end{align*} 
Since the second term becomes small as $\varepsilon \rightarrow 0$ independent of $n$, we see that this converges to $0$ indeed. 
\end{proof}
\end{proposition}

\begin{proposition}{
\label{convexcombi}The convex combinations $\sum_{j=1}^m \frac{1}{m} \text{ev}_{\lambda_j}$, where $\lambda_j$ are the solutions of $\lambda^m = 1$, are weak $^*$-limits in $\mathcal{S}(C(S^1))$ of sequences $(R_{n(m+1)}^*\tau_n)_{n\in \mathbb{N}}$ with $\tau_n \in \mathcal{P}(C(S^1)^{(n(m+1))})$.}
\begin{proof}
Take the polynomial $Q_\xi := \sum_{k} \xi_k z^{m-k} = \frac{1}{\sqrt{2}}( 1-z^{m})$, i.e. $\xi = \frac{1}{\sqrt{2}}(-1, 0, ..., 0, 1)$. Then the function \[g_m(t) := \abs{Q_\xi(e^{it})}^2 = 1 - \cos(mt)\] defines a pure state $\tau_{g_m}$ on $C(S^1)^{(m+1)}$ in the manner of Proposition~\ref{purestates1}. Likewise, due to Proposition~\ref{Pure state form}, the function $\frac{(g_m)^n}{\norm{(g_m)^n}_1}$ defines a pure state $\tau_n:=\tau_{\frac{(g_m)^n}{\norm{(g_m)^n}_1}}$ on $C(S^1)^{(n(m+1))}$.

Note that $g_m$ reaches its maximum in the $m$ points $\lambda_j$, denote the roots in between by $\mu_j$ (to be precise, each $\mu_j$ is the rotation of $\lambda_j$ by $\pi/m$). By symmetry of $g_m$, $\norm{(g_m\chi_{[\mu_j, \mu_{j+1}]})^n}_1 = \frac{1}{m} \norm{(g_m)^n}_1$. Hence \[\frac{(g_m)^n}{\norm{(g_m)^n}_1} = \frac{1}{m} \sum_{j=1}^m \frac{(g_m\chi_{[\mu_j, \mu_{j+1}]})^n}{\norm{(g_m\chi_{[\mu_j, \mu_{j+1}]})^n}_1},\]
and by applying Proposition~\ref{Bumpfunctions} on all these terms we conclude that \[R_{n(m+1)}^*\tau_n \xrightarrow{w^*} \sum_{j=1}^m \frac{1}{m} \text{ev}_{\lambda_j}. \qedhere\]
\end{proof}
\end{proposition}

Therefore, convex combinations of this type can indeed be approximated by pure states of the truncated circle. We will now use one more trick, which is to multiply the convergent sequence of pure states of the proposition above with the Radon-Nikodym derivative of another pure state, which results in a new sequence of pure states that converges to what we need.

\begin{proposition}{\label{productconv}Let $\lambda_1, ..., \lambda_m$ be the solutions of $\lambda^m = 1$, and take any $\varphi \in \mathcal{P}(C(S^1)^{(k)})$ such that $\frac{d\varphi}{d\lambda}(\lambda_j) \not=0$ for at least one $j$. Then \[ \sum_{j=1}^m \frac{\frac{d\varphi}{d\lambda}(\lambda_j)}{\sum_{i=1}^m\frac{d\varphi}{d\lambda}(\lambda_i) } \text{ev}_{\lambda_j}\] is the weak $^*$-limit of a sequence $(R_{k+n(m+1)}^*\chi_n)_{n \in \mathbb{N}}$ with $\chi_n \in \mathcal{P}(C(S^1)^{(k+n(m+1))})$.}
\begin{proof}
Consider the space of linear functionals $C(S^1)^*$. Observe that on $C(S^1)^*$:\begin{enumerate}
    \item The map \begin{align*}
    M_f^*: C(S^1)^* &\rightarrow C(S^1)^*\\
    \tau &\mapsto \tau \circ M_f,
\end{align*} where $M_f$ indicates multiplication by $f$, is weak$^*$-continuous. This is trivial, since if $\tau_n \xrightarrow{w^*} \tau$, then by definition $\tau_n(fg) \rightarrow \tau(fg)$ for all $g \in C(S^1)$ so $M_f^* \tau_n \xrightarrow{w^*} M_f^* \tau$.
\item If $\tau_n \xrightarrow{w^*} \tau$, then by definition also $\tau_n(1) \rightarrow \tau(1)$. As scalar multiplication is weak$^*$-continuous \[\frac{\tau_n}{\tau_n(1)} \xrightarrow{w^*} \frac{\tau}{\tau(1)},\] provided that these scalars are nonzero.
\end{enumerate}

Take $\varphi \in \mathcal{P}(C(S^1)^{(k)})$. For this proof, we will ease some notation by denoting the linear functional $g \mapsto \int_{S^1} fg \ \! d\lambda$ on $C(S^1)$ simply by $f$. As seen in the proof of Proposition~\ref{convexcombi}, if we define $g_m(t) = 1-\cos(mt)$ have that \[\frac{(g_m)^n}{\norm{(g_m)^n}_1} \xrightarrow{w^*} \sum_{j=1}^m \frac{1}{m} \text{ev}_{\lambda_j}.\]
If we apply observation 1 on this sequence with $M_{\frac{d\varphi}{d\lambda}}^*$, we get that \[\frac{(g_m)^n \frac{d\varphi}{d\lambda}}{\norm{(g_m)^n}_1} \xrightarrow{w^*}
\sum_{j=1}^m \frac{\frac{d\varphi}{d\lambda}(\lambda_j)}{m}\text{ev}_{\lambda_j}.\]
When $\frac{d\varphi}{d\lambda}(\lambda_j) \not=0$ for at least one $j$, all these are nonzero positive linear functionals on $C(S^1)$ so evaluating these functionals at $1$ gives a nonzero scalar. 

By observation 2, we can therefore conclude that \[\frac{(g_m)^n \frac{d\varphi}{d\lambda}}{\norm{(g_m)^n \frac{d\varphi}{d\lambda}}_1} \xrightarrow{w^*} \sum_{j=1}^m \frac{\frac{d\varphi}{d\lambda}(\lambda_j)}{\sum_{i=1}^m \frac{d\varphi}{d\lambda}(\lambda_i) } \text{ev}_{\lambda_j}.\]

Finally, the functional $f \mapsto \int_{S^1} f \frac{(g_m)^n \frac{d\varphi}{d\lambda}}{\norm{(g_m)^n \frac{d\varphi}{d\lambda}}_1} \ \! d\lambda$ is exactly $R_{k+(n(m+1))}^* \chi_n$ for $\chi_n$ the pure state on $C(S^1)^{(k+(n(m+1)))}$, defined via Proposition~\ref{Pure state form} by \[\frac{d\chi_n}{d\lambda} = \frac{(g_m)^n \frac{d\varphi}{d\lambda}}{\norm{(g_m)^n \frac{d\varphi}{d\lambda}}_1}. \qedhere \]
\end{proof}
\end{proposition}

We have completed all the steps that were described in the beginning of this subsection. Combined, that gives following proposition.

\begin{proposition}\label{convergence}
Given any state $\psi \in \mathcal{S}(C(S^1))$ and $\varepsilon >0$, there exists $N \in \mathbb{N}$ such that for any $n \geq N$ we can find a pure state $\chi_n \in \mathcal{P}(C(S^1)^{(n)})$ with $d(\psi, R_n^* (\chi_n)) < \varepsilon$, where $d$ is the Monge-Kantorovich metric on $\mathcal{S}(C(S^1))$. 
\begin{proof}
The proof of this proposition is nothing but the execution of the steps as described at the start of this section. There is one subtlety, however, so we do not omit the proof altogether.

Take $\psi \in \mathcal{S}(C(S^1))$ and $\varepsilon >0$. Using Lemma~\ref{approx} and the observation that the set  \[\left\{ \sum_{i=1}^n t_i \text{ev}_{\lambda_i} \ : \ 0 \leq t_i \leq 1, \sum_{i=1}^n t_i = 1, \lambda_i^n = 1 \right\}\] is dense in $\mathcal{S}(C(S^1))$ with respect to the weak *-topology, we can choose a (non-pure) state on $C(S^1)$ \[\rho = \sum_{j=1}^m  \frac{\frac{d\varphi}{d\lambda}(\lambda_j)}{\sum_{i=1}^m \frac{d\varphi}{d\lambda}(\lambda_i)} \text{ev}_{\lambda_j}\] with $\varphi$ a pure state in $\mathcal{P}(C(S^1)^{(m+1)})$ such that $d(\rho, \psi) < \varepsilon$.

According to Lemma~\ref{productconv} we can find a sequence $(\chi_n)_{n\in \mathbb{N}}$ with $\chi_n \in \mathcal{P}(C(S^1)^{(m+n(m+1))})$ such that $R_{m+n(m+1)}^*(\chi_n)$ converges to $\rho$ in the Monge-Kantorovich metric. However, we have to `fill in the gaps' to show that the distance of $\psi$ to the intermediate pure state spaces also shrinks arbitrarily small.

We can exploit Proposition~\ref{approx}, which shows that we can also find a pure state of higher `degree' $\varphi_l \in \mathcal{P}(C(S^1)^{(m+1+l)})$ such that $d(\rho_l, \psi) < \varepsilon$, for any $l\in \mathbb{Z}_{\geq 0}$. The argument above can then be repeated to find a sequence $\chi_{l,n} \in \mathcal{P}(C(S^1)^{(l+m+n(m+1))})$ converging to $\rho_l$. Choosing $l=0, ..., m$ results in a finite number of interlacing, separate sequences, which we can combine into one sequence of pure states $(\chi_n)_{n\in \mathbb{N}}$ such that $\chi_n \in \mathcal{P}(C(S^1)^{(n)})$. For this Frankenstein sequence, $R_n^* (\chi_n)$ might not converge, but there does exist an $N$ such that $n\geq N$ implies $d(R_n^*(\chi_n), \psi) < 2\varepsilon$ which proves the proposition.
\end{proof}
\end{proposition}

\subsection{Uniformity}
The result of the previous section means we can approximate all elements in $\mathcal{S}(C(S^1))$ by pullbacks of elements in $\mathcal{P}(C(S^1)^{(n)})$. In order to show Gromov-Hausdorff convergence, it remains to be shown that this approximation can be done uniformly so that $\mathcal{P}(C(S^1)^{(n)})$ forms an $\varepsilon$-net in $\mathcal{S}(C(S^1))$. A simple argument suffices, since $\mathcal{S}(C(S^1))$ is weak$^*$ compact.

\begin{proposition}
\label{epsilonnet}
For every $\varepsilon>0$, there exists $N\in \mathbb{N}$ such that for $n\geq N$, $R_n^* (\mathcal{P}(C(S^1)^{(n)}))$ forms an $\varepsilon$-net in $\mathcal{S}(C(S^1))$.
\begin{proof}
By the Banach-Alaoglu Theorem~\cite[Theorem~V.3.1]{Conway}, the unit ball of $C(S^1)^*$ is weak$^*$-compact. The set $\mathcal{S}(C(S^1))$, as a weak$^*$-closed subset of the unit ball, is then compact as well. Hence, given $\varepsilon > 0$, we can find a finite number of $\psi_1, ..., \psi_m \in \mathcal{S}(C(S^1))$ such that the balls $B_{\varepsilon}(\psi_i)$ cover $\mathcal{S}(C(S^1))$.

According to Proposition~\ref{convergence}, for each $\psi_i$ we can find $N_i \in \mathbb{N}$ such that for $ n \geq N_i$ there exists a pure state $\varphi_i \in \mathcal{P}(C(S^1)^{(n)})$ such that $d(\psi_i, R_n^*(\varphi_i)) < \varepsilon$. Thus we have that for any $\chi \in B_{\varepsilon}(\psi_i)$ and $n \geq N_i$, \[\operatorname{dist}(\chi, R_n^* (\mathcal{P}(C(S^1)^{(n)}))) \leq d(\chi, \psi_i) +\operatorname{dist}(\psi_i, R_n^* (\mathcal{P}(C(S^1)^{(n)}))) < \varepsilon + \varepsilon = 2 \varepsilon.\]

Now take any $\chi \in \mathcal{S}(C(S^1))$. Because \[\mathcal{S}(C(S^1)) \subseteq \bigcup_{i=1}^m B_{\varepsilon}(\psi_i),\]  $\chi$ must be an element of the ball $B_{\varepsilon}(\psi_i)$ for some $i$. For $n\geq N:= \max_i N_i$ we have, by the calculation above, that \[\operatorname{dist}(\chi, R_n^* (\mathcal{P}(C(S^1)^{(n)}))) < \varepsilon. \]
Hence $R_n^* (\mathcal{P}(C(S^1)^{(n)}))$ forms an $\varepsilon$-net in $\mathcal{S}(C(S^1))$ for $n \geq N$.
\end{proof}
\end{proposition}

\begin{corollary}{The set of measures on $S^1$ whose Radon-Nikodym derivatives with respect to the normalised Haar-measure on $S^1$ is of the form $c \prod_{j=1}^n (1-\cos(t-\theta_j))$ is dense in the set of Borel measures on $S^1$ with respect to the vague topology.}
\begin{proof}
Denote the set of positive Borel measures on $S^1$ by $\mathcal{M}^+(S^1)$. By Proposition~\ref{epsilonnet}, we see that for the set \[A:= \bigcup_{n=1}^\infty \left\{\mu \in \mathcal{M}^+(S^1) : \frac{d\mu}{d\lambda} =  c \prod_{j=1}^n (1-\cos(t-\theta_j)), c\in \mathbb{R}_{\geq 0}, \theta_1, \dots, \theta_n \in [0,2\pi)\right\},\] we have for any $\mu \in\mathcal{M}^+(S^1)$ that $\operatorname{dist}(A, \mu) = 0$. Hence $A$ is dense in $\mathcal{M}^+(S^1)$ with respect to the weak$^*$-topology, i.e. the vague topology.
\end{proof}
\end{corollary}

Proposition~\ref{epsilonnet} also concludes the proof of the Gromov-Hausdorff convergence of $\mathcal{P}(C(S^1)^{(n)})$ to \\ $\mathcal{S}(C(S^1))$ as metric spaces.

\begin{theorem}{\label{MainResult2}The metric spaces $\mathcal{P}(C(S^1)^{(n)})$ converge to the metric space $\mathcal{S}(C(S^1))$ in \\ Gromov-Hausdorff convergence.}
\begin{proof}
Corollary~\ref{epsilonisometry}, Lemma~\ref{R_n} and Proposition~\ref{epsilonnet} combined immediately give the result.
\end{proof}
\end{theorem}

\section{Recovering $S^1$}
It may come as a surprise that the limit of the spaces $\mathcal{P}(C(S^1)^{(n)})$ is $\mathcal{S}(C(S^1))$, and not $\mathcal{P}(C(S^1)) \cong S^1$. However, since by this identification $S^1$ is a subset of $\mathcal{S}(C(S^1))$, it must also be possible to recover $S^1$ as a Gromov-Hausdorff limit of a sequence of subsets of pure states. In this section we will demonstrate this. 

\begin{lemma}{\label{fejerkernels}The Fej\'er kernel rotated by $\lambda = e^{i\theta}$ \[f_n^\lambda(x)  = \sum_{\abs{k}\leq n-1}\left(1 - \frac{\abs{k}}{n} \right)e^{ik(\theta-x)},\] defines a pure state on $C(S^1)^{(n)}$ in the sense of Corollary~\ref{purestates2}}.\label{Fejer}\begin{proof}
Take the polynomial $\sum_{k=0}^{n-1} z^{k}$. Observe that \[(z-1)\sum_{k=0}^{n-1} z^{k} = z^{n} -1,\] and hence the roots of $\sum_{k=0}^{n-1} z^{k}$ are precisely the $n$th roots of unity with the exception of 1 itself. In fact, $\sum_{k=0}^{n-1} z^{k} = Q_\xi(z)$ where $\xi$ is the vector $\xi = \frac{1}{\sqrt{n}} (1, ..., 1) \in \mathbb{C}^n$, and $\abs{Q_\xi}^2$ defines a pure state. A simple calculation gives that \[\abs{Q_\xi(z)}^2 = \sum_{\abs{k}\leq n-1} \left(1-\frac{\abs{k}}{n} \right)z^k,\] hence the function $\sum_{\abs{k}\leq n-1}\left(1 - \frac{\abs{k}}{n}  \right)e^{ikx} = \sum_{\abs{k}\leq n-1}\left(1 - \frac{\abs{k}}{n}  \right)e^{-ikx} \in C(S^1)$ defines a pure state on $C(S^1)^{(n)}$. Rotations of this pure state are then also pure states by rotational invariance of $S^1$. \qedhere

\end{proof}
\end{lemma}

Compare these pure states $\tau_{f_n^\lambda}$ to the states on $C(S^1)$ denoted $\Psi_{x,N}^\sharp$ in~\cite[Section~5.4]{DAndrea}. The relation between these is that $R_n^*(\tau_{f_n^\lambda}) = \Psi_{\lambda, n}^\sharp$. These Fej\'er states recover the entire circle, which we will show in the next proposition. Note the similarity with~\cite[Proposition~5.11]{DAndrea} and~\cite[Proposition~5.12]{DAndrea}. The difference is that in this thesis we are talking about the \textit{intrinsic} distance on the truncated circle, so we have to add a a small step to move between the intrinsic distance and the distance on the whole spectral triple.

\begin{proposition}{\label{Localised states} Define the subsets $\mathcal{F}_n \subset \mathcal{P}(C(S^1)^{(n)})$ by \[\mathcal{F}_n := \{\tau_{f^\lambda_n} : \lambda \in S^1\},\] where $\tau_{f_n^\lambda}$ are states defined by Fej\'er kernels like in Lemma~\ref{fejerkernels}. Then the sequence of metric spaces $(\mathcal{F}_n, d_n)$ converges to $(S^1, d)$ in the Gromov-Hausdorff sense.}
\begin{proof}
Define the sets \[\mathfrak{R}_n = \{(\tau_{f^\lambda_n}, \lambda): \lambda \in S^ 1\} \subset \mathcal{F}_n \times S^1.\] Because the elements of $\mathcal{F}_n$ are labeled by $S^1$, the projections of $\mathfrak{R}_n$ onto the first and second coordinate are both surjective, making these sets total onto correspondences (Definition~\ref{correspondence}). We now want to show that the distortion of these correspondences converges to zero, in order to use Theorem~\ref{GHCorrespondence}. 

Because the Fej\'er kernel is a good kernel~\cite[Chapter~2]{hoffman2007banach}, the states $R_n^*(\tau_{f_n^\lambda})$ converge to $\text{ev}_\lambda$ in $\mathcal{S}(C(S^1))$ as $n \rightarrow \infty$. It thus follows immediately that \[\lim_{n\rightarrow \infty} d(R_n^*(\tau_{f_n^\lambda}), R_n^*(\tau_{f_n^\mu})) = d(\lambda, \mu),\] and because of Lemma~\ref{R_n} also \[ \lim_{n \rightarrow \infty} d_n(\tau_{f_{n}^\lambda}, \tau_{f_{n}^\mu}) = d(\lambda,\mu).\] 

We now want to estimate \[\sup_{\lambda, \mu \in S^1}\abs{d(R_n^*(\tau_{f_n^\lambda}), R_n^*(\tau_{f_n^\mu})) - d(\lambda,\mu)}.\] By definition, \[d(R_n^*(\tau_{f_n^\lambda}), R_n^*(\tau_{f_n^\mu})) = \sup_{g \in C^\infty(S^1)}\{ \abs{R_n^*(\tau_{f_n^\lambda})(g) - R_n^*(\tau_{f_n^\lambda})(g)} : \norm{g'}_{\infty} \leq 1\}.\] As $R_n^*(\tau_{f_n^\lambda})$ and $R_n^*(\tau_{f_n^\mu})$ are both states, we can freely subtract a constant function from $g$, so we might as well impose the extra condition \[d(R_n^*(\tau_{f_n^\lambda}), R_n^*(\tau_{f_n^\mu})) = \sup_{g\in C^\infty(S^1)}\{ \abs{R_n^*(\tau_{f_n^\lambda})(g) - R_n^*(\tau_{f_n^\lambda})(g)} : \norm{g'}_{\infty} \leq 1, \norm{g}_\infty \leq \pi\}.\]

Observe that for any smooth function $g$ with $\norm{g'}_\infty \leq 1$, $\norm{g}_\infty \leq \pi$, \begin{align*}
    \abs{R_n^*(\tau_{f_n^\lambda})(g) - R_n^*(\tau_{f_n^\mu})(g)} & \leq \abs{R_n^*(\tau_{f_n^\lambda})(g) - g(\lambda)} + \abs{R_n^*(\tau_{f_n^\mu})(g)-g(\mu)} + \abs{g(\lambda) - g(\mu)}\\
    & \leq \abs{R_n^*(\tau_{f_n^\lambda})(g) - g(\lambda)} + \abs{R_n^*(\tau_{f_n^\mu})(g)-g(\mu)} + d(\lambda, \mu),
\end{align*}
and furthermore if we denote the $\varepsilon$-neighbourhood of $\lambda$ in $S^1$ by $U^\varepsilon$,
\begin{align*}
     \abs{R_n^*(\tau_{f_n^\lambda})(g) - g(\lambda)} & \leq \frac{1}{2\pi} \int_{S^1} f_n^\lambda(x) \abs{g(x) - g(\lambda)} dx\\
     & =  \frac{1}{2\pi}\int_{S^1 \setminus U^\varepsilon} f_n^\lambda(x) \abs{g(x) - g(\lambda)} dx +  \frac{1}{2\pi}\int_{U^\varepsilon} f_n^\lambda(x) \abs{g(x) - g(\lambda)} dx\\
     & \leq \int_{S^1 \setminus U^\varepsilon} f_n^\lambda (x) dx + \varepsilon.
\end{align*}
Note that this last estimate is independent of $g$ and even of $\lambda$ by rotational invariance. Collecting the statements above, taking any $\varepsilon > 0$ gives \[\sup_{\lambda, \mu \in S^1}\abs{d(R_n^*(\tau_{f_n^\lambda}), R_n^*(\tau_{f_n^\mu})) - d(\lambda,\mu)} \leq  2 \int_{S^1 \setminus U^\varepsilon} f_n^\lambda (x) dx + 2 \varepsilon.\]
By the properties of the Fej\'er kernel~\cite[Chapter~2]{hoffman2007banach}, we may conclude that $\int_{S^1 \setminus U^\varepsilon} f_n^\lambda (x) dx$ converges to zero and therefore \[\lim_{n \rightarrow \infty} \operatorname{dist } \mathfrak{R}_n = 0, \] so by Theorem~\ref{GHCorrespondence} \[ \lim_{n \rightarrow \infty}d_{GH}(\mathcal{F}_n, S^1) = 0.\qedhere\]
\end{proof}
\end{proposition}

As a final remark, this proposition is comparable to a result proven by L. Glaser and A. Stern~\cite{GlaserStern}, which asserts that the pure state space of any spectral triple is the Gromov-Hausdorff limit of `localised' (not necessarily pure) states on the truncated spectral triple.

\textbf{Declaration.} \quad The author states that there is no conflict of interest.


\begin{thebibliography}{spmpsci}

\bibitem{burago}
D.~Burago, I.U.D. Burago, and S.~Ivanov.
\newblock {\em A Course in Metric Geometry}.
\newblock Crm Proceedings \& Lecture Notes. American Mathematical Soc., 2001.

\bibitem{caratheodory}
C.~Carath{\'e}odory and L.~Fej{\'e}r.
\newblock {\"U}ber den {Zusammenhang} der {Extremen} von harmonischen
  {Funktionen} mit ihren {Koeffizienten} und {\"u}ber den
  {Picard-Landau}’schen {Satz}.
\newblock {\em Rendiconti del Circolo Matematico di Palermo (1884-1940)},
  32(1):218--239, 1911.

\bibitem{Connes1994}
A.~Connes.
\newblock {\em Noncommutative Geometry}.
\newblock Elsevier Science, 1994.

\bibitem{ConnesSuijlekom}
A.~Connes and W.~D. {van Suijlekom}.
\newblock Spectral truncations in noncommutative geometry and operator systems.
\newblock {\em Communications in Mathematical Physics}, 2020.

\bibitem{Conway}
J.~B. Conway.
\newblock {\em A Course in Functional Analysis}.
\newblock Graduate Texts in Mathematics. Springer New York, 1994.

\bibitem{DAndrea}
F.~D’Andrea, F.~Lizzi, and P.~Martinetti.
\newblock Spectral geometry with a cut-off: topological and metric aspects.
\newblock {\em Journal of Geometry and Physics}, 82:18--45, 2014.

\bibitem{GlaserStern}
L.~Glaser and A.~B. Stern.
\newblock Reconstructing manifolds from truncations of spectral triples.
\newblock {\em Journal of Geometry and Physics}, 159:103921, 2021.

\bibitem{Scriptie}
E.~Hekkelman.
\newblock Truncated geometry on the circle.
\newblock Master's thesis, Radboud University Nijmegen, July 2021.

\bibitem{hoffman2007banach}
K.~Hoffman.
\newblock {\em Banach spaces of analytic functions}.
\newblock Courier Corporation, 2007.

\bibitem{Murphy}
G.~J. Murphy.
\newblock {\em C*-Algebras and Operator Theory}.
\newblock Elsevier Science, 1990.

\bibitem{rieffel1999metrics}
M.~A. Rieffel.
\newblock Metrics on state spaces.
\newblock {\em Documenta Mathematica}, 4:559--600, 1999.

\bibitem{Rieffel20041}
M.~A. Rieffel.
\newblock Gromov–hausdorff distance for quantum metric spaces.
\newblock {\em Mem. Am. Math. Soc.}, 168(796):1--65, 2004.

\bibitem{PapaRudin}
W.~Rudin.
\newblock {\em Real and Complex Analysis}.
\newblock Higher Mathematics Series. McGraw-Hill Education, 1987.

\bibitem{vansuijlekom2020gromovhausdorff}
W.~D. {van Suijlekom}.
\newblock {Gromov–Hausdorff} convergence of state spaces for spectral
  truncations.
\newblock {\em Journal of Geometry and Physics}, 162:104075, 2021.

\bibitem{vandermonde}
Z.~Yang, L.~Xie, and P.~Stoica.
\newblock Vandermonde decomposition of multilevel {Toeplitz} matrices with
  application to multidimensional super-resolution.
\newblock {\em IEEE Transactions on Information Theory}, 62(6):3685--3701,
  2016.

\end{thebibliography}
\end{document}